\documentclass[smallextended]{article}       % onecolumn (second format)
\usepackage{graphicx}
\usepackage{amssymb}
\usepackage{xypic}
\usepackage[english]{babel}
\usepackage{hyperref}

\newtheorem{theorem}{Theorem}[section]
\newtheorem{corollary}{Corollary}
\newtheorem{lemma}[theorem]{Lemma}
\newtheorem{proposition}{Proposition}

\newtheorem{remark}{Remark}

\newcommand{\fiz}{ { \textstyle{ \frac{1}{2} } } }

\newcommand{\bJ}{{\bf J}}

\newcommand{\bz}{{\bf z}}

\newcommand{\bx}{{\bf x}}

\newtheorem{condition}{Condition}%[section]

\begin{document}
\title{Symplectic maps: from generating functions to Liouvillian forms
\thanks{Funding: The author was supported by a grant from the \emph{Fondation du Coll\`ege de France}
and \emph{Total} under the research convention PU14150472, as well as the ERC Advanced Grant 
WAVETOMO, RCN 99285. 
Conflict of Interest: The author declares that he has no conflict of interest.}
}

\author{Hugo Jim\'enez-P\'erez
}

\maketitle

\begin{abstract}
In this article we introduce a new method for constructing implicit symplectic maps 
using \emph{special symplectic manifolds} and \emph{Liouvillian forms}. This method 
extends, in a natural way, the method of generating functions to 
1-forms which are globally defined on the symplectic manifold. 
The maps constructed by this method, are related to the symplectic Cayley's transformation 
and belong to a continuous space of dimension n(2n+1).
Applying the implicit map to the discrete Hamilton equations we obtain the 
\emph{generalized symplectic Euler scheme}.
We show the relations of the elements of this family with other discrete 
symplectic mapping, in particular 1) with the mappings obtained by generating 
functions of type I, II, and III and IV; 2)
with the symplectic Euler methods A and B; and 3) with the mid-point rule.
Moreover, we show the corresponding symplectic diffeomorphisms and their
Liouvillian forms on the product symplectic manifold.
We illustrate the details of the method 
in constructing two different families of implicit symplectic maps for $n=1$.
This is a geometrical method which overcomes the  
difficulties of the 
Hamilton-Jacobi theory and generating functions. 

%\keywords{Symplectic mappings, generating functions, Liouvillian forms, 
%symplectic integrators.}
%\subclass{Primary: 37M15, 65P10; Secondary: 53D22.}

\end{abstract}

\section{Introduction}

Symplectic maps can be constructed by using generating functions, and they were 
introduced by Poincar\'e  when looking for periodic orbits 
of second genus \cite{Poi99}. 
From the numerical point of view, 
symplectic maps are used for simulating Hamiltonian dynamics, based on the 
fact that the Hamiltonian flow is a one parameter subgroup of 
symplectic diffeomorphisms. In the classical construction 
of symplectic maps, the use of Darboux's coordinates in 
different stages of the procedure hides several interesting properties,
which arise with the discretization in time of Hamiltonian flows.
Other properties do not appear in the continuous case due to a classical
procedure of reduction from symplectic to contact geometry, 
which takes the Hamiltonian vector field into the Reeb field
on the reduced manifold \cite{LM87,Lib00}. 

In this paper, we address the problem of constructing implicit symplectic maps
with a different approach. Instead of 
using Hamilton-Jacobi theory to obtain a generating function,
we take advantage of the geometrical properties of both the Hamiltonian flow 
and the product symplectic manifold.  However, we avoid the use of Darboux's 
theorem and the classical evolutive Hamilton-Jacobi equation.
The reason is that the coordinates of the original problem are locally conjugated 
to the canonical coordinates by a symplectomorphism which is hidden by Darboux's 
theorem. In consequence, the classical construction integrates a modified problem 
in normal form using action-angle coordinates. 
Instead of Darboux's theorem we use \emph{special symplectic manifolds} \cite{Tul76,Tul77,ST72} 
and their \emph{Liouvillian forms} in order to track the state of the original coordinates 
at every stage of the procedure. 

This approach differs from the standard 
procedure of generating functions in two ways.
First, we consider Lagrangian submanifolds as the integral submanifolds of 
a prescribed Liouvillian form $\theta$ \cite{LM87,Lib00}, \emph{i.e.} submanifolds defined by the 
differential ideal $\mathcal I$ generated by $\theta$
\cite{BGG03}. 
Second, we use the complex structure on the product manifold to obtain 
the equations of the tangent space to the Lagrangian submanifold defined 
implicitly by the Liouvillian form. The advantage in the symplectic case is 
that the vertical bundle corresponds to the characteristic bundle to 
the Lagrangian submanifold. 
In our method, the Liouvillian form substitutes the generating function and the 
multiplication with the complex structure on the product manifold 
substitutes the role of the resolution of the Hamilton-Jacobi equation.
An additional projection of the Lagrangian submanifold offers the possibility to  
construct symplectic maps by a simple transversality condition.
For the linear case, this condition is weakened by the classical 
condition stating that the symplectic map must converge to the identity map 
on the diagonal of the product manifold.
All this procedure is constructed on the framework of special symplectic manifolds. 

The latter  condition is translated into a condition relating two different Liouvillian forms on 
the original symplectic manifold. Surprisingly, using this relation, the framework 
of special symplectic manifolds is no more 
necessary. In this way, we can construct symplectic maps just using 
Liouvillian forms, leading to the \emph{method of Liouvillian forms}.
This procedure reveals some remarkable properties and gives a different 
perspective on Liouvillian forms and generating functions. 
We give an interpretation 
of this method as a discrete version of the Cayley's transformation of a 
Hamiltonian matrix, leading naturally to a symplectic map. Finally, we construct a 
numerical scheme, which is the natural generalization of the Euler symplectic 
schemes. 

We summarize the rest of the paper. 
In Section 2 we give the main definitions and classical results on symplectic manifolds and the
product manifold that are necessary for our exposition. Special symplectic manifolds and Liouvillian forms are 
exposed in Section 3. 
In this section we introduce the generating functions of type I, II, III
and IV and their related Liouvillian forms expressed as elements of some special symplectic manifolds. 
We return systematically to these examples in the rest of the paper.
In Section 4 we expose the way we can construct symplectic maps 
using Liouvillian forms. 
Section 5 gives the structure of Liouvillian forms on the symplectic manifold
and the product symplectic manifold. In Section 6 we construct the implicit map. 
Section 7 describes all four steps of the 
method of Liouvillian forms and we state our main result. 
In addition, we relate these maps to the symplectic Cayley's transformation.
In Section 8 we apply the  implicit maps to Hamiltonian dynamics and we introduce the \emph{generalized implicit 
symplectic Euler scheme}. We show that the mid-point rule, and the symplectic Euler $A$ and $B$ maps are
particular cases of our scheme. Finally, in Section 9 we develop two examples 
of continuous families of implicit symplectic integrators: the first one using special symplectic manifolds 
and the second one applying the generalized symplectic Euler scheme.

\section{Symplectic manifolds and symplectic mappings}
In this section we recall some classical results and definitions in order to uniformize 
the notation. The results are stated in a geometrical fashion in preparation of the next section 
where we will set the framework of the method of Liouvillian forms. 
We assume the reader is familiar with the terminology of differential geometry 
and vector bundles. For an introduction the reader is referred to \cite{AM78,LM87,MR91}.

A \emph{symplectic manifold} is a $2n$-dimensional manifold $M$ equipped with 
a non-degenerated, skew-symmetric, closed 2-form $\omega$, such that
at every point $m\in M$, the tangent space 
$T_mM$ has the structure of a symplectic vector space. In addition, we say that
$(M,\omega)$ is an \emph{exact symplectic manifold} if there exists a 1-form 
$\theta$ such that $\omega=d\theta$; $\theta$ is called a Liouvillian form.
In what follows, all the symplectic manifolds will be exact.

An \emph{almost complex stucture} $J$ on a manifold $M$ is a smooth tangent
bundle isomorphism $J:TM\to TM$ covering the identity map on $M$ such
that for each point $z\in M$, the map $Jz=J(z): T_zM\to T_zM$ is a complex structure
on the vector space $T_zM$, it means an endomorphism on $T_zM$ such that 
$J^2(v)=J\circ J(v)=-v$ for every $v\in T_zM$.
We write $J^2=-I$ for simplicity.
A symplectic manifold with an almost complex structure possesses a Riemannian structure
$g$ which enables the definition of a symmetric positive definite form where %such that 
$\omega(\cdot,\cdot)= g(\cdot,J\cdot)$. We fix the 
Riemannian structure $g(\cdot,\cdot) = \langle\cdot,\cdot\rangle$ such that at every point on any symplectic 
manifold $x\in M$, the tangent space $(T_xM,\langle\cdot,\cdot\rangle_x)$ is
an inner product space. 

A submanifold $\Lambda\subset \tilde M$ is called \emph{Lagrangian} if the symplectic form restricted
to $\Lambda$ vanish $\omega|_\Lambda\equiv 0$.
The tangent bundle to the symplectic manifold $TM$ splits into the tangent $T\Lambda$ and the vertical $V\Lambda$
bundles of the Lagrangian submanifold in the form $TM=T\Lambda \oplus V\Lambda$, 
where $V\Lambda:=(T\Lambda)^\perp$ is the complementary subbundle to $T\Lambda$.
The complex structure $J$ and the non degeneracy of $\omega$ on $M$, let us rewrite
the vertical bundle as $V_x\Lambda = JT_x\Lambda$ for every $x\in \Lambda$.
We write this as
\begin{eqnarray*}
       T_x  M= T_x\Lambda \oplus J(T_x\Lambda).
\end{eqnarray*}

Let $(M,\omega)$ be a symplectic manifold. A diffeomorphism $f:M\to M$  
is said \emph{symplectic} or a \emph{symplectomorphism}
if $f^*\omega= \omega$, where $f^*$ denotes the pullback of 
the 2-form defined by
\begin{eqnarray}
   (f^*\omega)_z(u,v) = \omega_{f(z)}(Tf(u),Tf(v)),\quad z, f(z)\in M, u,v\in T_zM.
   \label{eqn:omega}
\end{eqnarray}
In this expression $Tf:TM\to TM$ is the tangent map of $f:M\to M$.
The set of all the symplectic diffeomorphisms on a symplectic manifold %$(M,\omega)$
form a group denoted by $Sp(M,\omega)$ and called the symplectic group
of $(M,\omega)$.

In order to construct symplectic maps we follow the classical construction. Define 
the product manifold of two 
copies of an exact symplectic manifold $(M,\omega)$ at times $t=0$ and $t=h$, which we denote by
$(M_1,\omega_1)$ and $(M_2,\omega_2)$, respectively. 
The manifold $\tilde M=M_1\times M_2$ with the canonical projections 
$\pi_i:\tilde M\to M_i$ for $i=1,2$ let us define the forms $\theta_\ominus$ and $\omega_{\ominus}$ on $\tilde M$ by 
\begin{eqnarray}
    \theta_{\ominus} = \pi_1^*\theta_1 - \pi_2^*\theta_2,\qquad
    \omega_{\ominus} = \pi_1^*\omega_1 - \pi_2^*\omega_2
    \label{eqn:def:sym}
\end{eqnarray}
We have the following facts (see \cite[sec 5.2]{AM78} for the proofs):
\begin{itemize}
    \item $(\tilde M,\omega_{\ominus})$ is a symplectic manifold of dimension $4n$.
    \item For any symplectic map $\phi:M_1\to M_2$, the graph of $\phi$, denoted by $\Lambda_\phi$,
        and defined by 
        \begin{eqnarray*}
            \Lambda_\phi = \left\{ \left( z,\phi(z) \right)\in \tilde M\ |\ z\in M_1,\phi(z)\in M_2 \right\}
        \end{eqnarray*}
        is a Lagrangian submanifold of $\tilde M$. It means $\omega_{\ominus}|_{\Lambda_\phi}\equiv 0$
    \item Since $\omega_{\ominus}=d\theta_{\ominus}$, then $\theta_\ominus$ is a locally closed form on 
        $\Lambda_{\phi}$. Applying Poincare's lemma,
        $\theta_{\ominus}$ is locally exact in a neighborhood of  $\Lambda_\phi$. Consequently, 
        there exists a function $S$ defined on
        the Lagrangian submanifold $\Lambda_{\phi}$ such that its differential concides with the
        restriction of the 1-form $\theta_{\ominus}$ to $\Lambda_{\phi}$.
        It means that 
        \begin{eqnarray}
	    dS|_{\Lambda_\phi} = \theta_\ominus|_{\Lambda_\phi}
        \end{eqnarray}
	    The function $S:\Lambda_{\phi}\to\mathbb R$ is called a \emph{generating 
	    function} for the symplectic map $\phi$;
    \item There exists and induced endomorphism on $T_x\tilde M$ which becomes the
        associated complex structure on $\omega_\ominus$ given by 
        \begin{eqnarray*}
            \tilde J = (T\pi_1)^T J_1 - (T\pi_2)^T J_2
        \end{eqnarray*}
        where $J_i$, are the associated complex structures to $\omega_i$, $i=1,2$.
\end{itemize}

Symplectic maps, generating functions, Liouvillian forms and Lagrangian submanifolds 
are closely related. For instance,  in a generic symplectic manifold $M$, any Lagrangian 
submanifold $\Lambda\subset M$ which is transverse to 
the fibers of the projection $\pi_{M}:T^* M\to M$ can be parameterized by a suitable 
(local) function $S$.
In contrast, Liouvillian forms are globally defined
and they do not need a particular parameterization, however topological properties 
of the underlying manifold like characteristic classes prevent the use of these 
forms in particular applications.

We rephrase  some classical properties of Lagrangian and symplectic submanifolds in 
term of the symplectic product manifold $\tilde{M}$.

\begin{lemma}
    Let $\Lambda\subset \tilde M$ be a Lagrangian submanifold and $\Phi\in Sp(\tilde{M},\omega_\ominus)$
    a symplectomorphism. We have the following facts:
    \begin{enumerate}
        \item The image of the Lagrangian submanifold under $\Phi$ is 
        again a Lagrangian submanifold of $\tilde{M}$.
    \item The projection $\pi_i(\Lambda)\subset M_i$ is a Lagrangian
    submanifold in $M_i$, $i=1,2$.
    \end{enumerate}
    \label{lem:1}
\end{lemma}
{\it Proof.} {\it 1)} Let $\bar\Lambda=\Phi(\Lambda)$ be the image of $\Lambda$ under the symplectomorphism $\Phi$.
  For all $y\in\bar\Lambda$ and $\xi,\eta\in T_{y}\bar\Lambda$ there exist $x\in\Lambda$ and $u,v\in T_x\Lambda$
  such that $y=\Phi(x)$, $\xi=T\Phi(u)$ and $\eta=T\Phi(v)$. Since $\Phi^*\omega_\ominus=\omega_\ominus$ we have 
  \begin{eqnarray*}
      (\omega_\ominus)_y(\xi,\eta)    = (\omega_\ominus)_x(u,v) = 0
  \end{eqnarray*}
  and $\bar\Lambda$ is a Lagrangian submanifold in $(\tilde M,\omega_\ominus)$.

{\it 2)} Since $\tilde M=M_1\times M_2$ is a product manifold, its tangent bundle splits
  into $T\tilde M= TM_1\oplus TM_2$. The subblundle $T\Lambda\subset T\tilde M$ is given under the 
  splitting by $T\Lambda=T\Lambda_1\oplus T\Lambda_2$. We have $T\pi_i(T\Lambda)= T\Lambda_i\subset TM_i$,   $i=1,2$.
  Due to the spltting $T\Lambda=T\Lambda_1\oplus T\Lambda_2$ then $\omega_\ominus|_{\Lambda}\equiv 0$ if
  and only if $\omega_i|_{\pi_i(\Lambda)}\equiv 0$ and consequently $\pi_i(\Lambda)\subset M_i$  are  
  Lagrangian submanifolds for $i=1,2$.
$\hfill\square$

\begin{lemma}
    Let $\phi_i\in Sp(M_i,\omega_i)$, $i=1,2$ be two symplectomorphisms. The induced
    diffeomorphism on $\tilde{M}$ by diagonal action on $M_1\times M_2$
    defined by $\Phi(\tilde{M})=\phi_1(M_1)\times\phi_2(M_2)$, is a symplectic
    diffeomorphism in $\tilde{M}$.
    \label{lem:3}
\end{lemma}
{\it Proof.} It is enough to show that $\Phi^*\omega_\ominus=\omega_\ominus$. By definition of the 
symplectic form and the fact that $\phi_i$, $i=1,2$ are symplectomorphisms we have 
\begin{eqnarray*}
      \omega_\ominus = \pi_1^*\omega_1 - \pi_2^*\omega_2  = \pi_1^*\circ \phi_1^*\omega_1 - \pi_2^*\circ \phi_2^*\omega_2 
      = (\phi_1\circ\pi_1)^*\omega_1 - (\phi_2\circ\pi_2)^*\omega_2.
\end{eqnarray*}
Applying the properties of the pull-back of the composition we obtain 
\begin{eqnarray*}
      \omega_\ominus  
      = (\phi_1\circ\pi_1)^*\omega_1 - (\phi_2\circ\pi_2)^*\omega_2 
       = \Phi^*\omega_\ominus 
\end{eqnarray*}
which gives the result.
$\hfill\square$

This last result gives us the possibility to consider Lagrangian submanifolds which do
not look like the graph of a diffeomorphism in $\tilde M$. It 
lets us work with Lagrangian submanifolds in mixed coordinates, or in other words
with generic implicit symplectic mappings. 

\section{Special symplectic manifolds and Liouvillian forms}

An important class of symplectic manifolds are the cotangent bundle to Riemannian
manifolds. In particular, they model the phase
space of many Hamiltonian mechanical systems, where the base coordinates correspond
to positions or configurations, and vertical coordinates correspond to momenta.

Let $\mathcal Q$ be a $C^\infty$ manifold. Consider the cotangent and tangent bundles 
$T^*\mathcal Q$ and $T\mathcal Q$, and its canonical projections
on the base manifold $\mathcal Q$ denoted by $\tau_{\mathcal Q}:T\mathcal Q\to \mathcal Q$ 
and $\pi_\mathcal Q:T^*\mathcal Q\to \mathcal Q$ respectively. 

The projection $\pi_\mathcal{Q}$ determines a natural map between the cotangent
bundle and the double cotangent bundle by its pullback $\pi_{\mathcal{Q}}^*:T^*\mathcal Q\to T^*(T^*\mathcal Q)$,
which sends 1-forms on $\mathcal Q$ to 1-forms on $T^*\mathcal Q$  by $\theta =\pi^*\eta$ 
for every 1-form $\eta$ defined on $\mathcal Q$. Formally $\eta$ is a section 
of the cotangent bundle that we denote by $\eta\in \Gamma(T^*\mathcal Q)$, 
where $\Gamma(T^*\mathcal Q)$ denotes the 
space of smooth sections on the cotangent bundle.\footnote{We use the notation 
$\beta\in_\Gamma E$ instead of  $\beta\in\Gamma(E)$ for simplifying the notation in the diagrams. } 
\begin{eqnarray*}
    \xymatrix{
        &   \theta\in_\Gamma\ T^*(T^*\mathcal Q)  \ar[dr]^{\pi_{T^*\mathcal Q}}
                        & \\
        \eta\in_\Gamma\ T^*\mathcal Q \ar[ur]^{\pi^*_\mathcal Q} \ar[dr]^{\pi_\mathcal Q} \ar[rr]^{id} & &  T^*\mathcal Q \ar[dl]_{\pi_\mathcal Q}\ \ni_\Gamma\ \eta\\
        & \mathcal Q &  
    }
\end{eqnarray*}

The form $\theta=\theta_\eta$, induced by the identity morphism $\theta_\eta=(\pi^*_\mathcal Q)_\eta(\eta)$ 
for every $\eta \in \Gamma(T^*\mathcal Q)$, is called 
the \emph{Liouville form} on $T^*\mathcal Q$ and is alternatively defined by its action on the tangent bundle
by the equation
  \begin{eqnarray*}
    \langle v, \theta_\eta\rangle_{T^*\mathcal Q} &=& \langle T\pi_\mathcal Q(v), \eta  \rangle_{\mathcal Q}, 
        \qquad v\in T_\eta T^*\mathcal Q,\quad \eta\in \Gamma(T^*\mathcal Q), \quad \theta_{\eta}\in T^*_\eta(T^*\mathcal Q).
  \end{eqnarray*}

\begin{eqnarray*}
    \xymatrix{
        &   v\in TT^*\mathcal Q \ar[dl]_{T\pi_\mathcal Q} \ar[dr]^{\tau_{T^*\mathcal Q}}
                       \ar@{--}[rr]|{ \langle \cdot ,\cdot \rangle_{T^*\mathcal Q}} & & T^*(T^*\mathcal Q) \ar[dl]_{\pi_{T^*\mathcal Q}}  \ni \theta_{\eta}  \\
        T\mathcal Q \ar[dr]^{\tau_\mathcal Q} \ar@{--}[rr]|{ \langle \cdot ,\cdot \rangle_{\mathcal Q}} & &  T^*\mathcal Q \ar[dl]_{\pi_\mathcal Q}\ni_\Gamma\eta&\\
        & \mathcal Q &  &
    }
\end{eqnarray*}

\begin{remark}
The double cotangent bundle $(T^*(T^*\mathcal Q), T^*\mathcal Q, \pi_{T^*\mathcal Q})$
is defined by the projection $\pi_{T^*\mathcal Q}:T^*(T^*\mathcal Q)\to T^*\mathcal Q$. Consequently, 
the canonical Liouville form $\theta\in\Gamma(T^*(T^*\mathcal Q))$ is a section in $T^*(T^*\mathcal Q)$ 
corresponding to the inclusion of sections from $\Gamma(T^*\mathcal Q)$ into $\Gamma(T^*(T^*\mathcal Q))$.
This inclusion is also interpreted as the identity map. The fact that the Liouville form is
a section of the double cotangent bundle is not evident when we work
with symplectic vector spaces and it has been the source of many misunderstandings 
in some areas like the numerical analysis.
\end{remark}

When we shall have the occasion to deal with cotangent bundles of  different manifolds, we will denote 
the Liouville form on $T^*\mathcal N$ by $\theta_\mathcal N$. 

The cotangent bundle $T^*\mathcal Q$ inherits a natural symplectic structure $\omega=d\theta_\mathcal Q$ 
such that the couple $(T^*\mathcal Q,\omega)$ becomes a symplectic manifold.
Unfortunately, many geometrical properties of symplectic manifolds lost significance in mechanical systems
due to a missing physical interpretation. Tulczyjew proposed in \cite{Tul76} the study of symplectic manifolds which
are diffeomorphic to cotangent bundles by means of \emph{special symplectic manifolds} (see also \cite{Tul77,ST72}).

    A \emph{special symplectic manifold} is a quintuple $( M, \mathcal Q, \theta, \pi, \varphi )$
    where $\pi: M \to \mathcal Q$ is a fibre bundle, $\theta$ is a 1-form on $M$ and
    $\varphi:M \to T^*\mathcal Q$ is a symplectic diffeomorphism such that $\pi=\pi_\mathcal Q\circ \varphi$,
    and $\theta=\varphi ^*\theta_\mathcal Q$. 
\begin{eqnarray*}
    \xymatrix{
        & \theta\in_\Gamma T^*M\ \ar[dl]_{\pi_{M}} &  
             & \ar[ll]_{\varphi^*}  T^*(T^*\mathcal Q) \ar[dl]_{\pi_{T^*\mathcal Q}} \ni_\Gamma \theta_{\mathcal Q}  \\
        M\ \ar[dr]_{\pi }  \ar[rr]^{ \varphi } & & \ T^*\mathcal Q \ar[dl]_{\pi_\mathcal Q} &   \\
       &  \mathcal Q & & 
    }
\end{eqnarray*}

The pair $(M,d\theta)$ is a symplectic manifold. We call the 1-form $\theta=\varphi^*\theta_\mathcal Q$
a \emph{Liouvillian form}\footnote{Libermann and Marle \cite{LM87}
and Libermann \cite{Lib00} define Liouvillian forms in the more general context
of bundle morphisms. 
We use a simplified definition for working with special symplectic manifolds.} on $M$. More generally 
we say that a 1-form $\eta$ on an even dimensional symplectic  manifold $(M,\omega)$
is \emph{Liouvillian} or of \emph{Liouvillian type} if $d\eta$ is a symplectic form on $M$, or equivalently if $(M,d\eta)$
is a symplectic manifold.

Let $K\subset \mathcal Q$ be a submanifold and $S:K\to \mathbb R$ a function on $K$. The \emph{Lagrangian 
submanifold generated} by $S$ in the manifold $({M},d\theta)$ is defined by the equation
$\langle v, \theta\rangle = \langle T\pi(v), dS \rangle$ where $v\in T{M}$ and $\tau_{{M}}(v)=m$,
in the following way 
\begin{eqnarray}
    \Lambda = \left\{ m\in{M} | \pi(m)\in K,\langle v, \theta\rangle =\langle T\pi(v),dS\rangle \right\}.
    \label{eqn:lag}
\end{eqnarray}

We suppose $K\subset \mathcal Q$ is an embedded submanifold and we write it as an inclusion $i:K \to\mathcal Q$. 
In this case $i^*(dS)\in \Gamma(T^*K)$ and $dS\in\Gamma(T^*\mathcal Q)$. If in addition, we
suppose that $i:K\to\mathcal Q$ is an immersion, $\pi$ defines a submersion and 
the submanifold $\Lambda$ is well defined. The following 
graph shows this situation
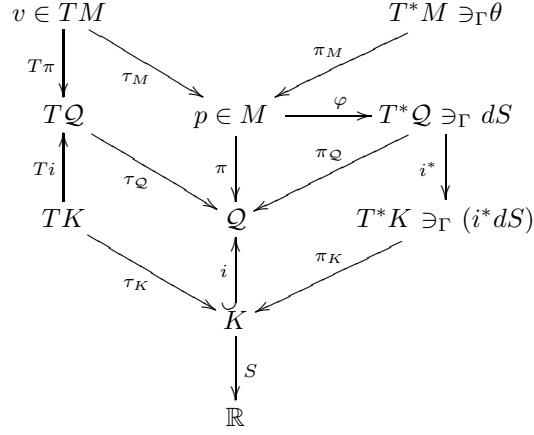
\begin{eqnarray*}
    \xymatrix{
        v\in T  M\ \ar[dr]_{\tau_{  M}} \ar[d]_{T\pi} & &  T^*  M\ \ar[dl]_{\pi_{  M}} {\ni_\Gamma} \theta \\
       T\mathcal Q \ar[dr]_{\tau_{\mathcal Q}}  & p\in   M\ \ar[d]_{\pi }  \ar[r]^{ \varphi } & 
       T^*\mathcal Q\ {\ni_\Gamma}\ dS \ar[dl]_{\pi_{\mathcal Q}} \ar[d]_{i^*} \\
       TK \ar[u]^{Ti} \ar[dr]_{\tau_{K}} & \mathcal Q &    T^*K\ {\ni_\Gamma}\ (i^*dS) \ar[dl]_{\pi_{K}}  \\
      & K  \ar[d]^{S} \ar@{^{(}->}^{i}[u] & \\
      &\mathbb R&
    }
\end{eqnarray*}

Any special symplectic manifold $(M,\mathcal Q,\theta,\pi,\varphi)$ determines a natural splitting 
of the tangent bundle $TM=T\Lambda\oplus V\Lambda$ using the Liouvillian form $\theta=\varphi^*\theta_{\mathcal Q}$ 
and the projection $\pi=\pi_\mathcal Q\circ\varphi$. The projection
keeps track of the deformation of the tangent (horizontal) bundle, and the Liouvillian form 
recovers the deformation of the vertical bundle.

The following result relates the symplectic manifolds of interest in this work: $(\tilde{M},\omega_\ominus)$ 
and $(T^*(\mathcal Q_1\times\mathcal Q_2), d\theta_{\mathcal Q_1\times\mathcal Q_2})$, where 
$\theta_{\mathcal Q_1\times\mathcal Q_2}$ is the Liouville form on $\mathcal Q_1\times\mathcal Q_2$.
\begin{proposition}
    Define coordinates $(q,p,Q,P)\in \tilde{M}$ in the product manifold such that
    $(q,p)\in M_1$ and $(Q,P)\in M_2$. 
    Let $E_1:\tilde{M}\to T^*(\mathcal Q_1\times\mathcal Q_2)$ be the linear map
    given by 
    \begin{eqnarray}
        E_1:\tilde{M} &\to& T^*(\mathcal Q_1\times\mathcal Q_2) \nonumber\\
          (q,p,Q,P) & \mapsto & (q,-Q,p,P).
        \label{eqn:e1}
    \end{eqnarray}
    Then $E_1$ is a symplectomorphism.
\end{proposition}
{\it Proof.} 
Define coordinates 
$(x,X,y,Y)\in T^*(\mathcal Q_1\times \mathcal Q_2)$ where $(x,X)\in\mathcal Q_1\times \mathcal Q_2$.
We have $(x,X,y,Y)=(q,-Q,p,P)$ and the Liouville form in these coordinates becomes 
$\theta_{\mathcal Q_1\times \mathcal Q_2}= ydx + YdX$. Taking the differential we have
$$E_1^*(dy\wedge dx + dY\wedge dX) = dp\wedge dq - dP\wedge dQ=\omega_\ominus$$
as we want to show.

$\hfill\square$

$E_1$ is known as the \emph{canonical symplectomorphism} between the
product manifold $(\tilde{M},\omega_\ominus)$ and the 
contangent bundle $(T^*(\mathcal Q_1\times\mathcal Q_2), d\theta_{\mathcal Q_1\times\mathcal Q_2})$. 

\begin{corollary}
\label{cor:1}
  Let $\Phi\in Sp(\tilde{M},\omega_\ominus)$ be a symplectomorphism on $\tilde{M}$ and consider 
  $E_1$ as below. Then the diffeomorphism $\Psi:\tilde{M}\to 
  T^*(\mathcal Q_1\times\mathcal Q_2)$ given by 
  \begin{eqnarray}
  \Psi=E_1\circ\Phi
  \label{eqn:symp}
  \end{eqnarray}
  is symplectic.
\end{corollary}
{\it Proof.} Direct from the fact that  $\Psi: \tilde M\to T^*(\mathcal Q_1\times \mathcal Q_2)$
is the composition of two symplectomorphisms 
\begin{eqnarray*}
    \xymatrix{
        \tilde M  \ar[dr]_{\Psi}  \ar[r]^{\Phi} & \tilde M \ar[d]^{E_1}\\
        & T^*(\mathcal Q_1\times\mathcal Q_2).
    }
\end{eqnarray*}
$\hfill\square$

The quintuple $(\tilde{M},\mathcal Q_1\times\mathcal Q_2, \theta,\pi,\Psi )$,
where $\Psi$ is defined in (\ref{eqn:symp})
becomes a special symplectic manifold on $\mathcal Q_1\times\mathcal Q_2$.

\begin{eqnarray*}
    \xymatrix{
        & \ T^*\tilde M\ \ar[dl]_{\pi_{\tilde M}} &  & \ar[ll]_{\Psi^*}  T^*T^*(\mathcal Q_1\times\mathcal Q_2) 
         \ar[dl]^{\quad\pi_{T^*(\mathcal Q_1\times\mathcal Q_2)}}  \\
        \ \tilde M\ \ar[dr]_{\pi }  \ar[rr]^{ \Psi } & & \ T^*(\mathcal Q_1\times\mathcal Q_2) \ar[dl]^{\quad\pi_{\mathcal Q_1\times\mathcal Q_2}} &   \\
      &  \mathcal Q_1\times\mathcal Q_2 & & 
    }
\end{eqnarray*}

We recall that any Lagrangian submanifold $\Lambda\in T^*(\mathcal Q_1\times\mathcal Q_2)$ defines 
a Lagrangian submanifold in $\tilde M$ by the symplectomorphism
$\Psi^{-1}(\Lambda)$. 

\subsection{Examples of special symplectic manifolds}
\label{sec:ex}
We consider some simple symplectomorphisms between the product of two symplectic manifolds and 
the cotanget to the product of two configuration spaces $\tilde M\to T^*(\mathcal Q_2\times \mathcal Q_1) $.
These symplectomorphisms define four special symplectic manifolds and
we will be interested in their Liouvillian forms. 
Using symplectic coordinates on the manifolds $(x_i,X_i,y_i,Y_i)\in T^*(\mathcal Q_1\times\mathcal Q_2)$ and 
$(q_i,p_i,Q_i,P_i)\in\tilde M$ we consider four symplectic diffeomorphisms 

\begin{enumerate}
    \item $E_1(q_i,p_i,Q_i,P_i)=(q_i,-Q_i,p_i,P_i)$ corresponding to the canonical 
    symplectic diffeomorphism from  $\tilde M\to T^*(\mathcal Q_2\times \mathcal Q_1)$.
    The Liouville form on $T^*(\mathcal Q_2\times \mathcal Q_1)$ corresponding to 
    $\theta_{\mathcal Q_1\times\mathcal Q_2}$ is pulled back to the Liouvillian 
    form $\theta=(E_1^*)\theta_{\mathcal Q_1\times\mathcal Q_2}$
    and it corresponds to the form $\theta_\ominus$ on $\tilde M$.  
    For instance, the Liouvillian form and the projection are  
    \begin{eqnarray*}
       \theta_{\ominus}=p_idq_i -P_idQ_i \qquad{\rm and} \qquad \pi (q_i,p_i,Q_i,P_i)=(q_i,-Q_i).
    \end{eqnarray*}
    \item $\Psi_{\imath\imath}(q_i,p_i,Q_i,P_i)=(q_i,P_i,p_i,Q_i)$. This diffeomorphism gives the
    Liouvillian form and projection given by 
    \begin{eqnarray*}
	\theta=p_idq_i +Q_idP_i\qquad{\rm and}\qquad \pi (q_i,p_i,Q_i,P_i)=(q_i,P_i).
    \end{eqnarray*}
    \item $\Psi_{\imath\imath\imath}(q_i,p_i,Q_i,P_i)=(Q_i,p_i,-P_i,-q_i)$. This diffeomorphism gives the
    following Liouvillian form and projection
    \begin{eqnarray*}
	\theta=-q_idp_i -P_idQ_i\qquad{\rm and}\qquad\pi (q_i,p_i,Q_i,P_i)=(Q_i,p_i).
    \end{eqnarray*}
    \item $\Psi_{\imath v}(q_i,p_i,Q_i,P_i)=(-p_i,P_i,q_i,Q_i)$. This diffeomorphism gives the
    following Liouvillian form and projection 
    \begin{eqnarray*}
	\theta=-q_idp_i +Q_idP_i,\qquad {\rm  and}\qquad \pi (q_i,p_i,Q_i,P_i)=(-p_i,P_i).
    \end{eqnarray*}
\end{enumerate}

Since the Liouvillian form on $\tilde M$ corresponds to a class of symplectomorphisms modulo 
symplectic rotations we need both elements, the Liouvillian form and the projection for fixing the 
symplectomorphism between these manifolds. 

The symplectomorphims $E_1,\Psi_{\imath\imath},\Psi_{\imath\imath\imath},\Psi_{\imath v}$ will be revisited 
in the next section where we related them to generating functions.

\section{Symplectic maps from Liouvillian forms}

The usual way to construct symplectic maps on the product manifold 
$(\tilde M,\omega_\ominus)$ is using 
generating functions $S:\Lambda\to\mathbb R$ defined on some 
Lagrangian submanifold $\Lambda\subset\tilde M$; it is an inverse problem.
This inverse problem is solved using Hamilton-Jacobi theory for estimating the
characteristic bundle which contains, as a subbundle, the vertical bundle to the 
Lagrangian submanifold. 
We deal with this problem in a more direct way using Liouvillian forms which define 
Lagrangian submanifolds as their integral submanifolds. 
The transformation between the vertical and the tangent bundles
is given by the complex structure associated to the symplectic form.
In this way we avoid the solution of both the Hamilton-Jacobi 
equation and the generating function. However, in both procedures, 
we must select what type of symplectic maps we are looking for.
In our case, we are interested on symplectic maps
adapted for constructing numerical schemes which differ from those 
used for studying periodic orbits in an essential way. The task in this
section is to characterize
Lagrangian submanifolds adapted for constructing numerical schemes.

\begin{remark}
  Symplectic maps for numerical schemes and those used for studying 
  periodic orbits  solve variational problems with different boundary conditions. 
  The former consider the minimization of the 
  action integral along paths joinning two different fixed points; the latter
  consider closed paths with prescribed period $T>0$. This implies that 
  Poincar\'e's differential form\footnote{It is the differential of the 
  so called Poincar\'e's generating function} introduced in \cite{Poi99}
  for studying bifurcations of periodic orbits is not well suited for 
  constructing numerical schemes. A detailed study of this fact is given in 
  \cite{Jim15f}.
\end{remark}

To obtain a symplectic map for constructing numerical schemes 
we need to recover a symplectic vector space % submanifold $N\subset\tilde M$ 
from the information encoded in the tangent space to the Lagrangian submanifold  
$\Lambda\subset \tilde M$.
In fact, we state a stronger condition. Let 
$(\tilde M,\mathcal Q_1\times\mathcal Q_2,\theta,\pi,\Psi)$
be the special symplectic manifold defined in the last section and
consider the image $N=\pi(\tilde M)$ of the projection $\pi=\pi_{\mathcal Q_1\times\mathcal Q_2}\circ \Psi$ 
as a submanifold by the adapted inclusion  $i:N{\hookrightarrow}\tilde M$.
If we can give a symplectic structure $\omega_N$ on $N$
such that $\omega_\ominus = (\pi^*)\omega_N$, and $\omega_N=(i^*)\omega_\ominus$ 
then the integral 
submanifolds of the Liouvillian form $\theta$ on $(\tilde M,\omega_\ominus)$
are adapted for the construction of symplectic integrators.

If the symplectic submanifold $N\subset \tilde M$ belongs to a symplectic path joining
$z(t)\in M_1\times M_1$ to $z(t+h)\in M_1\times M_2$, then the submanifold $(N,\omega_N)$
must converge to the original symplectic manifold
\begin{eqnarray}
    \lim_{h\to 0} (N,\omega_N) = (M_1,\omega_1)=(M_2,\omega_2).
\end{eqnarray}
We summarize these requirements in the following:
\begin{condition}
\label{cond:1}
    Consider the tangent space to the projected manifold at a point $\pi(x)\in N$, $x\in \tilde M$, 
    given by $T_{\pi(x)}N=T\pi(T_x\tilde M)$. Denote by $\bx\in \Delta\subset \tilde M$ 
    a point on the diagonal of the product manifold. If 
    \begin{enumerate}
     \item the restriction of $\omega_\ominus$ to the submanifold $N=\pi(\tilde M)$ is a symplectic form
     on $N$, and, 
    \item $\left.T\pi\left(T_x\tilde M\right)\right|_{x=\bx} = id$, 
    is the identity,
    \end{enumerate}
    the Liouvillian form $\theta$ of the special symplectic manifold
    $(\tilde M,\mathcal Q_1\times \mathcal Q_2, \theta,\pi,\Psi)$ defines
    a Lagrangian submanifold adapted for the construction of a symplectic scheme.
\end{condition}

Instead of recovering the symplectomorphism $\Psi$ which
defines the special symplectic manifold, we use the fact that 
Lagrangian submanifolds are the integral submanifolds
to Liouvillian forms. 
Let $\theta_1$ be a Liouvillian form on  the $2n$-dimensional symplectic 
manifold $(M,\omega)$. For every symplectomorphism $\phi: M\to M$, the 
form $\theta_2$ given by the pullback $\theta_1=\phi^* \theta_2$ is again 
Liouvillian on $(M,\omega)$. 
They are different Liouvillian forms producing
the same symplectic structure $\omega=d\theta_1=d\theta_2$ on $M$.
On the product manifold $(\tilde M,\omega_\ominus)$, the 
form $\theta = \pi_1^*\theta_1 - \pi_2^*\theta_2$ 
is Liouvillian on $\tilde M$ since $d\theta = \omega_\ominus$.

Using point $2)$ in Lemma \ref{lem:1} and Condition \ref{cond:1},
we will induce the Liouvillian form on the 
product manifold $(\tilde M,\omega_\ominus)$ by Liouvillian forms $\theta_1$ and $\theta_2$ on the 
original symplectic manifold $(M,\omega)$ defining complementary Lagrangian 
submanifolds. 
For this, we consider the complex structure $J$ associated to $\omega$.
There exists a (tautological) symplectomorphism $\bJ:M\to M$ such that 
its tangent map is exactly the complex structure $T\bJ = J:TM\to TM$. This symplectomorphism
is attached to every $\omega$.

\begin{theorem}
\label{teo:1}
For every exact symplectic manifold $(M,\omega)$, there exists at least a Liouvillian form $\theta_1$
on the original symplectic manifold $(M,\omega)$
such that the Liovillian form 
\begin{eqnarray}
  \theta = \pi_1^*\theta_1-\pi_2^*\theta_2,\qquad \theta_1=\bJ^*(\theta_2),
\end{eqnarray}
on the  product manifold $(\tilde M, \omega_\ominus)$ satisfies Condition 1.
\end{theorem}
We need additional elements for proving this theorem which follows from Lemma \ref{lem:7}
below. 
For the moment, we will relate
the Liouvillian forms of the examples given in section \ref{sec:ex} with the 
classical generating functions.

\begin{lemma}
\label{lem:II:III}
    The Liouvillian form of the symplectomorphisms $E_1, \Psi_{\imath\imath}, \Psi_{\imath\imath\imath}, \Psi_{\imath v}$ given in 
    section \ref{sec:ex} are associated to the generating functions of type I, II, III and IV, respectively.
\end{lemma}
{\it Proof.} We perform the same computations for every symplectomorphism using the 
Liouvillian form $\theta$ on the product manifold $(\tilde M,\omega_\ominus)$ 
and the projection $\pi = \pi_{(\mathcal Q_1\times\mathcal Q_2)}\circ\Psi$
for the corresponding $\Psi$.
\begin{enumerate}
    \item $E_1$ produces the Liouvillian form $\theta=p_idq_i -P_idQ_i$, which is locally 
    equivalent to the differential of a function $S:\tilde M\to\mathbb R$ with $S=S(q_i,Q_i)$ which 
    is a generating function of type I. It Defines 
    a Lagrangian submanifold in $\tilde M$  by 
    $$\Lambda=\{(\hat q_i,\hat p_i,\hat Q_i,\hat P_i) \in \tilde M| 
     \textstyle{ p_i=\frac{\partial S}{\partial q_i}, P_i=-\frac{\partial S}{\partial Q_i}} \}.$$
     It does not produces a map adapted for construct symplectic 
     integrators since the projection $\pi(q_i,p_i,Q_i,P_i)=(q_i,-Q_i)$ does not give symplectic coordinates.
    \item $\Psi_{\imath\imath}$ gives the
    Liouvillian form $\theta=p_idq_i +Q_idP_i$, which is locally 
    equivalent to the differential of a function $S=S(q_i,P_i)$, which is of type II. It defines 
    a Lagrangian submanifold in $\tilde M$  by
$$\Lambda_{\imath\imath}=\left\{(\hat q_i,\hat p_i,\hat Q_i,\hat P_i)\in \tilde M|\ \textstyle{p_i=\frac{\partial S}{\partial q_i}, 
Q_i=\frac{\partial S}{\partial P_i} }\right\}.$$
    This Liouvillian form is adapted for constructing a symplectic integrator since the projection 
    $\pi (q_i,p_i,Q_i,P_i)=(q_i,P_i)$ gives symplectic coordinates in the diagonal.
    \item $\Psi_{\imath\imath\imath}$ produces the
    Liouvillian form $\theta=-q_idp_i -P_idQ_i$ which is locally 
    equivalent to the differential of a function $S=S(Q_i,p_i)$, which is of type III. It defines 
    a Lagrangian submanifold in $\tilde M$  by
	$$\Lambda_{\imath\imath\imath}=\left\{(\hat q_i,\hat p_i,\hat Q_i,\hat P_i) \in \tilde M|\ 
\textstyle{ q_i=-\frac{\partial S}{\partial p_i}, 
P_i=-\frac{\partial S}{\partial Q_i}} \right\}.$$
    This Liouvillian form is adapted for constructing a symplectic integrator since the projection 
	$\pi (q_i,p_i,Q_i,P_i)=(Q_i,p_i)$ gives symplectic coordinates in the diagonal.
    \item $\Psi_{\imath v}$ gives the
    Liouvillian form $\theta=-q_idp_i +Q_idP_i$ which is locally 
    equivalent to the differential of a function $S=S(p_i,P_i)$, which is of type IV. It defines 
    a Lagrangian submanifold in $\tilde M$  by
	$$\Lambda_{\imath v}=\left\{(\hat q_i,\hat p_i,\hat Q_i,\hat P_i) \in \tilde M|\ 
\textstyle{ q_i=-\frac{\partial S}{\partial p_i}, 
Q_i=-\frac{\partial S}{\partial P_i}} \right\}.$$
    This Liouvillian form is not adapted for constructing a symplectic integrator since the projection 
    $\pi (q_i,p_i,Q_i,P_i)=(-p_i,P_i)$ does not give symplectic coordinates.
\end{enumerate}

\section{Structure of Liouvillian forms}
Generating functions for constructing symplectic maps are those of type \emph{II} and \emph{III}, 
which is not the case for generating functions of type \emph{I} and \emph{IV}. 
This happens in the same way for Liouvillian forms. In this section we study the structure
of Liouvillian forms for characterizing those adapted for constructing
implicit symplectic integrators. 
We recall some facts about the Liouvillian forms on a generic exact symplectic manifold 
$(M,\omega)$. 

\subsection{Liouvillian forms on an exact symplectic manifold}
Consider an exact symplectic manifold $(M,d\theta)$, where $\theta$ is 
a generic Liouvillian form. For every function  $F:M\to \mathbb R$ 
on $M$ the differential form  $\theta + dF$ is again a Liouvillian form
since $d(\theta + dF) = \omega$. It means,
the set of Liouvillian forms is invariant under the addition of 
exact 1-forms. This comes from the fact that $d^2(F)$ corresponds to 
a symmetric tensor
$S={\rm Hess}(F)$. 
This generic fact implies that the space of Liouvillian forms
is infinite-dimensional. 

For studying the structure of Liouvillian forms, we want to understand 
the structure of the symplectic form in a generic basis. 
Given local coordinates $z=\{z_i\}_{i=1}^{2n}\in M$, a basis of 
$T^*_zM \cong \Omega^1(M,\omega)$ is given by $dz=\{dz_i\}_{i=1}^{2n}$,
and a basis of $T_zM \cong \wedge^1(M,\omega)$ is given by 
$\partial_z=\left\{\frac{\partial}{\partial z_i}\right\}_{i=1}^{2n}$.
These are dual bases satisfying $dz_i\left( \frac{\partial}{\partial z_j}\right)= \frac{\partial z_i}{\partial z_j} = \delta_{ij}$.
In the basis $\{dz_i\}_1^{2n}$, the symplectic form $\omega$ is given by
\begin{eqnarray}
   \omega = \sum_{i,j}^{2n} \fiz J_{ij} dz_j\wedge dz_i, \qquad J_{ij}\in\mathbb R.
\end{eqnarray}
The matrix $J=\left[J_{ij}\right]$ satisfies $J^2 = -I_{2n}$, and it is the 
representation of the complex 
structure associated to $\omega$ in the basis $\{dz_i\}_1^{2n}$.
Remark that $J$ is an antisymmetric matrix $J+J^T=0_{2n}$.

A Liouvillian form $\theta$  on $(M,\omega)$ is expressed in the basis $\{dz_i\}_{i=1}^{2n}$
by $\theta = \sum_i \alpha_{i}(z)dz_i$ $z=(z_1,z_2,\cdots,z_{2n})\in M$,
where $\alpha_i:M\to\mathbb R$ are smooth functions of $z$. Since 
\begin{eqnarray*}
    \omega = d\theta = \sum_{i,j} \frac{\partial \alpha_j(z) }{\partial z_i} dz_i\wedge dz_j,\qquad z=(z_1,z_2,\cdots,z_{2n})\in M,
\end{eqnarray*}
and writing  $A(z)=\left[A_{ij}(z) \right] = \left[\frac{\partial \alpha_{i}(z)}{\partial z_j}\right]$ we have the decomposition 
in symmetric and antisymmetric components by 
\begin{eqnarray*}
     A_s(z)= \fiz \left( A(z) +A^T(z) \right) \qquad{\rm and} \qquad 
     A_a(z) = \fiz\left( A(z) - A^T(z)\right). 
\end{eqnarray*}
The condition for $\theta$ to be a Liouvillian form is $A_a(z)\equiv \fiz J$, and 
denoting $A_s(z) = S(z)$ to remark it is a symmetric matrix,
we have  $A(z) = S(z) + \fiz J$. 
Moreover,  the symmetric part $S(z)$ belongs to the kernel of the differential 
$\theta \stackrel{d}{\mapsto}\omega$.
We have proved the following:
\begin{lemma}
\label{lem:5}
   The set of Liouvillian forms on an exact symplectic 
   manifold $(M,\omega)$ 
   is given in local coordinates by  1-forms $\theta=\sum_i \alpha_i(z)dz_i$ 
   where the matrix $A(z)=\frac{\partial \alpha_j(z) }{\partial z_i}$ has the structure $A(z) = S(z) + \frac12 J$, where $S(z)= S^T(z)$
   and $J$ is the complex structure associated to $\omega$ in the local coordinates $\{z_i\}_1^{2n}$.
\end{lemma}

Every exact symplectic manifold $(M,\omega)$ possesses a 
natural Liouvillian form $\theta_0$ given by the complex structure associated 
to the symplectic form $\omega$. The Liouvillian form $\theta_0$ has no
symmetric component and we call it the \emph{elementary} or \emph{basic Liovillian form}. 
\begin{corollary}
    The elementary 
    Liouvillian form in Darboux's coordinates $z=(q,p)\in (M,\omega)$ is given by 
    \begin{eqnarray}
        \theta_0 = \fiz\sum_{i=1}^n\left(p_idq_i-q_idp_i\right).
    \end{eqnarray}
\end{corollary}

Consider the tautological symplectomorphism $\bJ$ as in Theorem \ref{teo:1}.
\begin{lemma}
\label{lem:7}
    The elementary Liouvillian form $\theta_0$ is invariant under the action of 
    the complex structure, {\it i.e.} $\theta_0=\bJ^*\theta_0$.
\end{lemma}
{\it Proof.} Given local coordinates $z\in( M,\omega)$ the elementary Liouvillian form 
is given by $\theta_0=\fiz dz J z$. Consequently\footnote{Remark that $dz$ is a covector 
and it transforms by the transpose $dz\mapsto dz J^T$.}
\begin{eqnarray}
    \bJ^*\theta_0 = \fiz (dzJ^T)J(Jz) = \fiz dz Jz = \theta_0. 
\end{eqnarray}
$\hfill\square$

This lemma gives a Liouvillian form on an exact symplectic manifold, proving 
Theorem \ref{teo:1}.
The Liouvillian form $\theta=\pi_1^*\theta_0 - \pi_2^*\theta_0$ induces a natural 
symplectomorphism which will be associated to the mid-point rule. We will show this fact in the 
last section where we will construct numerical integrators from the symplectic 
maps coming from Liouvillian forms.

For studying the way a Liouvillian form induces a symplectic integrator 
we restrict our study to the set
of Liouvillian forms with linear coefficients. It becomes a linear space over 
$\mathbb R$ with finite dimension. 
In the rest of this paper we replace the matrix $A(z)$ for a linear form 
$Az = (S+\fiz J)z$ where $A$ and $S$ are constant matrices in $\mathbb M_{2n\times 2n}(\mathbb R)$.

\begin{lemma}
    The space of Liouvillian forms with linear coefficients on a 2n-dimensional 
    exact symplectic manifold 
    $(M,\omega)$  has dimension $n(2n + 1)$. 
\end{lemma}
{\it Proof.} Applying Lemma \ref{lem:5}, we deduce that the dimension of the space of Liouvillian forms 
is exactly the dimension of the space $Sym(2n)$ of symmetric $2n\times 2n$ matrices $S=S^T$, given by 
$\dim Sym(2n)= \frac12(2n)(2n+1)=n(2n+1).$
$\hfill\square$

\begin{remark}
   The space of Liouvillian forms with linear coefficients on a symplectic manifold 
   $(M,\omega)$ is isomorphic to $\mathfrak{sp}(M,\omega)$ as an affine space. 
\end{remark}

\subsection{Liouvillian forms on the product manifold}
We are interested on  the space of Liouvillian forms with 
linear components on the product manifold $(\tilde M,\omega_\ominus)$ given by 
\begin{eqnarray}
   \theta = \pi_1^*\theta_1 - \pi_2^*\theta_2,\qquad \theta_i\in\Gamma(T^*M_i), i=1,2,
   \label{eqn:liouv}
\end{eqnarray}
where $\theta_i,i=1,2$ are represented in local coordinates by 
$\theta_i = dz(\fiz J+S_i)z$, $i=1,2$
(see Lemma \ref{lem:5}).

\begin{lemma}
   \label{lem:key1}
   Consider local coordinates $(z,Z)\in(\tilde M,\omega_\ominus)$, where $z\in M_1$ and $Z\in M_2$.
   Liouvillian forms on $(\tilde M,\omega_\ominus)$ given in (\ref{eqn:liouv}) have a 
   representation $\theta = (dz,dZ) \tilde A(z,Z)^T$ for the matrix $\tilde A = \tilde S + \frac12 \tilde J$, 
   where $\tilde A\in \mathbb M_{4n\times 4n}(\mathbb R)$
   is a symmetric matrix with the form 
\begin{eqnarray*}
    \tilde S = \left( 
     \begin{array}{c c}
       S_1 & 0_{2n} \\
       0_{2n} & -S_2
    \end{array}
    \right),\qquad S_1, S_2, 0_{2n}\in \mathbb M_{2n\times 2n}(\mathbb R), \quad S_i=S_i^T, i=1,2.
\end{eqnarray*}
\end{lemma}
{\it Proof.}
Consider a point on the product manifold $(z,Z)^T\in (\tilde M,\omega_\ominus)$ written 
in local coordinates of the factors
$z\in M_1$ and $Z\in M_2$. Using Lemma 
\ref{lem:5}, $\theta$ has a representation 
\begin{eqnarray}
\theta = (dz, dZ) \tilde A \left( \begin{array}{c} z\\ Z\end{array}\right), \qquad \tilde A=\tilde S + \fiz\tilde J,
\label{eqn:lem:5}
\end{eqnarray}
where $\tilde S$ is a symmetric matrix of size $4n\times 4n$.
The key property is the decomposition of the tangent bundle 
$T\tilde M=TM_1\oplus TM_2$, since the representation of Liouvillian 
forms of type (\ref{eqn:liouv}) is given by a matrix with two blocks  
\begin{eqnarray*}
    \tilde A = \left( 
     \begin{array}{c c}
       \fiz J + S_1  & 0_{2n} \\
       0_{2n} & -(\fiz J + S_2)
    \end{array}
    \right) =
    \left( 
     \begin{array}{c c}
       S_1  & 0_{2n} \\
       0_{2n} & -S_2
    \end{array}
    \right) + \fiz
    \left( 
     \begin{array}{c c}
       J & 0_{2n} \\
       0_{2n} & - J 
    \end{array}
    \right)= \tilde S +\fiz \tilde J
\end{eqnarray*}
proving the lemma.
$\hfill\square$

The (affine) space of Liouvillian forms with linear coefficients on $(\tilde M,\omega_\ominus)$ has dimension 
$\dim Sym(4n) = 2n(4n+1)$ and the dimension of the subset of forms given by 
(\ref{eqn:liouv}) 
is $2 \dim Sym(2n) = 2n(2n + 1)$.
However, just a subset of this space is adapted for the  construction of  numerical 
schemes.
We will be interested in the particular case $S=S_1=-S_2$ whose projection 
satisfies Condition \ref{cond:1}. 
We will return to this argument when 
we study Hamiltonian systems in Section \ref{sec:ham}.

\section{The implicit map from the projection of the Lagrangian submanifold}

The linearization of $\theta$ give us the local expresion of the 
vertical bundle\footnote{The null space or the kernel of de projection $T\pi$.}
$V\Lambda\subset T\tilde M$
to the Lagrangian submanifold $\Lambda$ and we want to recover the 
tangent bundle $T\Lambda \subset T\tilde M$.
Since the tangent bundle $T\tilde M$ accepts a decomposition by 
$T\tilde M = V\Lambda \oplus T\Lambda$ and the tangent bundle 
is mapped into the vertical bundle by $V_z\Lambda = \tilde J T_z\Lambda$
for every $z\in\Lambda$, then we recover locally the tangent spaces using the complex 
structure
$T_z\Lambda = \tilde J^{T} ( V_z\Lambda$). We project this subbundle with 
the tangent projection
\begin{eqnarray}
TN=T\pi(T\Lambda)=T\pi(\tilde J^T( V\Lambda)).  
\end{eqnarray}
\begin{lemma}
   Let $(\tilde M,\mathcal Q_1\times\mathcal Q_2,\theta,\pi,\Psi)$ be a special 
   symplectic manifold  and $\Lambda\subset \tilde M$ be an integral submanifold for
   $\theta = \pi_1^*\theta_1 - \pi_2^*\theta_2$ where $\theta_i = dz(\fiz J+S_i)z$. Then 
   for every $x\in\Lambda$, the subspaces $T_{\pi_i(x)}\Lambda_i:= T\pi_i (\tilde J^T(V_x\Lambda)) \subset T_xM_i$, $i=1,2$
   are Lagrangian with local expression $L_i=(\fiz I -JS_i)z$,
   satisfying $(\fiz J+S_i)z=0$, $i=1,2$.
   \label{lem:10}
\end{lemma}
{\it Proof.} 
Applying point \emph{2)} from Lemma \ref{lem:1}, the projection of the Lagrangian submanifold $\Lambda\in\tilde M$,
by $\pi_i(\Lambda)= \Lambda_i\subset M_i$ is a Lagrangian submanifold in $M_i$, $i=1,2$.
In local coordinates, the expressions 
\begin{eqnarray}
  L_i=\left(\fiz I-JS\right)z\qquad{\rm and}\qquad \left(\fiz J+S_i\right)z=0, \quad i=1,2,
\end{eqnarray}
are the equations of the tangent space to the 
integral submanifold $\Lambda_i\subset M_i$ defined by the Liouvillian form $\theta_i=dz(\fiz J+ S_i)z$.
They define the same tangent space since $\Lambda$ is an integral submanifold
for $\theta = \pi^*\theta_1 - \pi_2^*\theta_2$ and $\Lambda_i$ are integral submanifolds for
$\theta_i$.
$\hfill\square$

We define the linear space given by the sum $V = T_{\pi_1(x)}\Lambda_1 + T_{\pi_2(x)}\Lambda_2$
as vector subspaces of $T\tilde M$. Using Lemma \ref{lem:10} we can write $V=L_1 + L_2$.
We want that $V$ satisfies Condition \ref{cond:1},
however, point \emph{2)} in Condition \ref{cond:1} is enough for a symplectic map
as we will see below.

\begin{lemma}
Define the implicit map induced by the sum of the linear spaces $V=L_1 +L_2$
from Lemma \ref{lem:10}, given by
\begin{eqnarray}
   \rho(z,Z) = \left( \fiz I_{2n} - JS_1 \right)z + 
       \left( \fiz I_{2n} + JS_2 \right) Z.
\end{eqnarray}
Then, $\rho(z,z) = z$, if and only if $S_1=S_2$.
\end{lemma}
{\it Proof.} Evaluating in $(z,z)\in \Delta\subset \tilde M$ we have $\rho(z,z)=z +J(S_2-S_1)z$ 
which produces the identity on the diagonal, if and only if $J(S_2-S_1)z=0$.
$\hfill\square$

The implicit map induced by the projection of the lagrangian submanifold $\Lambda$
adapted for the  construction of symplectic maps is given in local coordinates of the product manifold 
$(z,Z)\in \tilde M$ by 
\begin{eqnarray}
   \rho(z,Z) = \fiz \left( z + Z\right) + b \left( Z-z\right), 
  \label{eqn:forms}
\end{eqnarray}
where $b\in \mathbb M_{2n\times 2n}(\mathbb R)$ is a Hamiltonian matrix $b^TJ+Jb=0$.

\section{The method of Liouvillian forms}
The construction of implicit symplectic maps using the method of Liouvillian
forms is obtained by restructuring the construction of the last section
avoiding the explicit use of the special symplectic manifold. We recover all the information
from two different (but related) Liouvillian forms and 
the complex structure.

{\it 1)} For any exact symplectic manifold $(M,\omega)$ 
consider two copies $(M_i,\omega_i)$, $i=1,2$ for the construction of the product manifold 
$(\tilde M,\omega_\ominus)$ and select a Liouvillian  form $\theta$ on $\tilde M$.

We define the Liouvillian form satisfying 
Condition \ref{cond:1} in the following way: 
fix a Liouvillian form $\theta_1$ with linear coefficients on $(M_1,\omega_1)$ and express it
in local coordinates by $\theta_1= dz(\fiz J+S)z$;
define a second Liouvillian form on $(M_2,\omega_2)$ by $\theta_2=dZ(\fiz J-S)Z$.
The form 
\begin{eqnarray}
   \theta = \pi_1^*\theta_1 - \pi^*_2\theta_2, 
\end{eqnarray}
is Liouvillian on the product manifold. 
The expression of the vertical bundle in local coordinates is  
\begin{eqnarray*}
  \left( 
     \begin{array}{c}
        z\\
        Z
     \end{array}
  \right)
  \mapsto
  \left( 
    \begin{array}{c c}
       \fiz J + S & 0_{2n}\\
       0_{2n} & -\fiz J  + S
     \end{array}
  \right)   
  \left( \begin{array}{c}
    z\\
    Z
    \end{array}
  \right), \qquad S\in M_{2n\times 2n}(\mathbb R), S=S^T.
\end{eqnarray*}

{\it 2)} 
The Pffafian equation $\theta=0$ defines implicitly and locally the Lagrangian 
submanifold $\Lambda$. We use the complex structure $\tilde J$ associated to the symplectic form $\omega_\ominus$ 
of the product manifold to obtain the equations of the tangent spaces by $T\Lambda = \tilde J^{T} V\Lambda$. 
In local coordinates %$(z,Z)\in (\tilde M,\omega_\ominus)$ 
we have 
%If we define $\hat z = $
\begin{eqnarray*}
  \left( 
     \begin{array}{c}
        z\\
        Z
     \end{array}
  \right)
  \mapsto
  \left( 
    \begin{array}{c c}
       \fiz I_{2n} - b & 0_{2n}\\
       0_{2n} & \fiz I_{2n}  + b
     \end{array}
  \right)   
  \left( \begin{array}{c}
    z\\
    Z
    \end{array}
  \right), \qquad b=JS\in M_{2n\times 2n}(\mathbb R), S=S^T.
\end{eqnarray*}

{\it 3)} Project this Lagrangian subspace on the tangent bundle of the original 
symplectic manifold $TM$ with $T\pi_1(T\Lambda) + T\pi_2(T\Lambda)$
as linear subspaces. 
We will have the induced implicit map 
\begin{eqnarray}
     \rho (z,Z)  &=& \left( \fiz I_{2n} - b \right)z + 
       \left( \fiz I_{2n} + b \right) Z,
       \label{eqn:proj:coord}
\end{eqnarray}

{\it 4)} Verify if $$TM = T\pi_1(T\Lambda) \oplus T\pi_2(T\Lambda)=T\Lambda_1\oplus T\Lambda_2,$$ 
which holds when the projected Lagrangian submanifolds are complementary. 
Since $z$ and $Z$ are two different close points,
we only check if the implicit map restricted to the diagonal is the identity $\rho|_{\Delta}=I_{2n}$. 
The map (\ref{eqn:proj:coord})
satisfies this condition by construction.

\begin{theorem}
\label{teo:2}
    Let $(M,\omega)$ be an exact symplectic manifold of dimension $2n$ and $U\subset M$ a start-shaped open
    subset containing the points $z,Z\in U$. Define an implicit map $\rho:U\times U\to U$ 
    by expression (\ref{eqn:proj:coord}). % and denote by $\bar z=\rho(z,Z)$ its image. Then
    The implicit map $(z,Z)\mapsto Z= \rho(z,Z)$ is symplectic.
\end{theorem}
{\it Proof.} We prove this result using implicit differentiation. Consider the implicit  mapping $\Psi$ given by 
\begin{eqnarray}
    \Psi(z,Z) = Z - \rho (z, Z) = 0.
    \label{eqn:th1}
\end{eqnarray}
Implicit differentiation of (\ref{eqn:th1}) gives 
\begin{eqnarray}
        \frac{\partial \Psi(z,Z)}{\partial z} = -\fiz I_{2n} + b\quad
    {\rm and }\quad  
       \frac{\partial \Psi(z,Z)}{\partial Z} = \fiz I_{2n} - b.
    \label{eqn:mats}
\end{eqnarray}
Denoting the partial derivatives of $\Psi$ by 
\begin{eqnarray*}
    B_1 =  -\frac{\partial \Psi(z,Z)}{\partial z}\quad
    {\rm and }\quad 
    B_2 = \frac{\partial \Psi(z,Z)}{\partial Z},
\end{eqnarray*}
the amplification matrix of the linearized system \cite{STW92}
is given by $B=B_2^{-1}\circ B_1$, and $\Psi$ is symplectic if the matrix 
$B$ is symplectic. 
We use the fact that the transpose  
$B^T=B_1^{T}\circ B_2^{-T}$ of a symplectic matrix is symplectic to tranform the 
symplecticity condition into 
$ B_2^{-1}\circ B_1 J B_1^{T}\circ B_2^{-T} = J$, or equivalently into
$ B_1JB_1^T - B_2JB_2^T = 0.$
This condition is satisfied since $B_1=B_2$, and the implicit map is symplectic 
as we want to prove.
$\hfill\square$

\begin{remark}
   The converse is in general not true since we are characterizing only symplectic maps solving a
   variational problem with particular boundary conditions.
\end{remark}

We can consider the point $\bar z=\rho(z,Z)$ as an intermediate point on the  path
joining $z$ to $Z$ minimizing some variational problem. To develop this idea, we introduce
some additional terminology. An implicit map $\phi:U\times U\to U$
is called \emph{consistent} if there exist two explicit maps $\psi_1,\psi_2:U\to U$ and a point $\bar z\in U$, such that 
\begin{eqnarray}
    \bar z = \psi_1(z)\qquad{\rm and}\qquad \bar z = \psi_2(Z).
\end{eqnarray}
We say that $\bar z$ is the \emph{point of consistency} and $\psi=\psi_2^{-1}\circ\psi_1:U\to U$
its \emph{consistency map}. It is an explicit well defined map. We say that $\phi$ \emph{interleaves 
a symplectic map} if its consistency map is symplectic. 

We need a result from the Weyl's extension of Cayley's transformation. The interested readers 
are referred to \cite{Wey46} for the proof. 
\begin{lemma}[Generalized Cayley's transformation]
  \label{lem:cay}
  If the non-exceptional matrices\footnote{A matrix $A\in GL(n)$ is said 
to be \emph{non-exceptional} if $\det(I+A)\neq0$, where $I$ is the 
identity matrix in $GL(n)$.} ${\bf H}$ and ${\bf S}$ are connected by the relations 
  \begin{eqnarray*}
      {\bf S} &=& (I-{\bf H})(I+{\bf H})^{-1} = (I+{\bf H})^{-1}(I-{\bf H})\\
      {\bf H} &=& (I-{\bf S})(I+{\bf S})^{-1} = (I+{\bf S})^{-1}(I-{\bf S})
  \end{eqnarray*}
and $G$ is any matrix, then ${\bf S}^TG{\bf S}=G$ if and only if ${\bf H}^TG + G{\bf H}= 0$. 
\end{lemma}

For the symplectic case, we fix the arbitrary matrix $G=J$ to the complex structure on $TM$.
Now we can relate the induced implicit map $\rho$ with Cayley's transformation in the following
\begin{proposition}
  With the same hypotheses as Theorem \ref{teo:2}, consider the explicit 
  map $\psi=\psi_2^{-1}\circ\psi_1:U\to U$ 
  as the consistency map
  associated to the implicit map $\rho:TU\times TU \to TU $ given in (\ref{eqn:proj:coord}).
  Then, the consistency map $\psi$ is symplectic and it corresponds to Cayley's transformation 
  of the matrix $2JS$.
\end{proposition}
{\it Proof.}
Consider $\rho(z, Z)$ as a 
linear combination of two explicit linear maps coming from $z$ and $Z$,
in the form $\rho(z, Z) = \frac12\left( T\psi_1(z) + T\psi_2(Z)\right)$. 
From expression (\ref{eqn:proj:coord}), we can write explicitely
\begin{eqnarray}
    T\psi_1(z) = (I_{2n} - 2JS)z \qquad{\rm and}\qquad T\psi_2(Z)= (I_{2n}+2JS)Z.
\end{eqnarray}
Since $S\in M_{2n\times 2n}(\mathbb R)$ is a symmetric matrix, the matrix ${\bf H}=2JS$
is Hamiltonian. In this case, the consistency map associated to $\rho(z,Z)$ is given by
$\psi=\psi_2^{-1}\circ \psi_1$, whose linearization gives 
\begin{eqnarray}
    T\psi = T\psi_2^{-1}\circ T\psi_1 =(I_{2n}+{\bf H})^{-1}(I_{2n} - {\bf H}).
\end{eqnarray}
This is the symplectic Cayley's transformation of the Hamiltonian matrix ${\bf H}$. Applying Weyl's 
Lemma we obtain immediately that the consistency map $\psi:U\ni z \mapsto Z=\psi(z)\in U$ is symplectic.
$\hfill\square$

\begin{corollary}
   For every symplectic mapping $\psi\in Sp(M,\omega)$ on an exact symplectic manifold,
   there exists a Liouvillian form $\theta_\psi$ adapted to $\psi$.
\end{corollary}
{\it Proof.} 
Suppose $\psi\in Sp(M,\omega)$ is represented at every point $m\in M$ by a 
non-exceptional matrix ${\rm\bf S}\in Sp(2n)$.
There exists an associated Hamiltonian matrix 
given by the symplectic Cayley's transformation
   ${\rm\bf H} = (I_{2n} - {\rm\bf S})(I_{2n} + {\rm\bf S})^{-1}$.
Consequently, the linearization of the Liouvillian form $\theta_\psi$ adapted 
to $\psi$ is given by some square 
matrix $R = \frac12 J\left( I_{2n} + {\bf H} \right)$
in the form $\overline{\theta_\psi} = dz R z$. Moreover, the implicit map 
$\rho:TU\times TU \to TU$ 
in coordinates $(z,Z)\in U\times U$ is defined 
in terms of ${\bf H}$ by 
\begin{eqnarray}
     \rho (z,Z)  &=& \left( \frac12 I_{2n} - \frac12{\bf H} \right)z + 
       \left( \frac12 I_{2n} + \frac12{\bf H} \right) Z.
       \label{eqn:proj:coord:H}
\end{eqnarray}
$\hfill\square$

Consider the matrix $G$ in Weyl's lemma \ref{lem:cay} as a bilinear form $u^T G v$.
The trivial case $G=I_{2n}$ gives the Euclidean metric and the Cayley's 
transformation relates symmetric with antisymmetric matrices. The case $G=J$
is the symplectic form and it relates symplectic with Hamiltonian matrices. 
In \cite{JVR17} we show that this property with an additional dependency on 
a parameter discretizing time produce reversible maps.

\section{Hamiltonian systems and implicit symplectic integrators}
\label{sec:ham}

The numerical simulation of Hamiltonian vector fields 
is made by numerical schemes called \emph{symplectic algorithms} or \emph{symplectic integrators}.
We introduce the notation and terminology for the application 
of Liouvillian forms in symplectic integration.

A \emph{Hamiltonian system} $(M,\omega,X_H)$ is a vector field
$X=X_H$ on a symplectic manifold $(M,\omega)$ such that inner product 
of the symplectic form and the vector field satisfies $i_{X_H}\omega = -dH$,
for a differentiable function $H:M\to\mathbb R$.

The set of solutions $\varphi^t_H:M\times \mathbb R\to M$  of the Hamiltonian vector field 
$X_H$ is called \emph{the flow of $X_H$} and it is defined by 
\begin{eqnarray}
   \frac{d}{dt} \varphi^t_H(z) = X_H(\varphi^t_H(z)).
\end{eqnarray}
The flow is a 1-parameter subgroup of symplectic diffeomorphisms with parameter $t\in\mathbb R$.
For every point $z_0\in M$, the solution $\varphi_H^t(z_0)$ is the integral curve of $X_H$ passing by $z_0$ at 
time $t=0$. 

A \emph{symplectic algorithm} with stepsize $h$, is the numerical 
approximation $\psi_h$ of the 
map $\varphi^t_H$, for $t=h$ fixed, which is smooth with respect to $h$ and $H$, 
and preserves the symplectic form $(\psi_h)^*\omega = \omega$. 
Consider an open set $U\in M$ and two points on the flow of $X_H$, say $z_t=\varphi^t(z_0)$
and $z_{t+h} = \varphi^{(t+h)}_H(z_0)$.
By the group property of the flow, it is enough to perform the analysis around $t=0$,
 then the points 
will be denoted by $z_0$ and $z_h$.
 \begin{theorem}
    Let $(M,\theta,X_H)$ be a Hamiltonian system on an exact symplectic manifold. 
    Consider a convex open set $U\subset M$ containing the points $z_0$ and $z_h$ on the
    flow of $X_H$. If an implicit map $\rho:U\times U\to U$ is defined by 
     \begin{eqnarray}
     \rho (z_0,z_h)  &=& \left( \fiz I_{2n} - b \right)z_0 + 
       \left( \fiz I_{2n} + b \right)z_h,
       \label{eqn:proj:coord:b}
 \end{eqnarray}
 where $b$ is a Hamiltonian matrix in $M_{2n\times 2n}(\mathbb R)$, 
    then, the map $$z_h = z_0 + h X_H\circ\rho(z_0,z_h)$$ 
    is symplectic. 
    \label{teo:main}
\end{theorem}
{\it Proof.}
Since $X_H$ is a Hamiltonian vector field, invariant under symplectic
transformations, it suffices that $\rho(z_0,z_h)$ be a symplectic map.  
Applying Theorem \ref{teo:2}, we obtain
the desired result. 
$\hfill\square$

An alternative proof of the theorem is using implicit differentiation. Consider the implicit  
mapping $\Psi$ given by 
\begin{eqnarray}
    \Psi(z_0,z_{h}) = z_h - z_0 - hX_H \circ \rho (z_0, z_h) = 0
    \label{eqn:th2}
\end{eqnarray}
as in the proof of Teorem \ref{teo:2}. 
Implicit differentiation of (\ref{eqn:th2}) using  the chain rule  gives 
\begin{eqnarray}
        B_1 := -\frac{\partial \Psi(z_0,z_h)}{\partial z_0} &=& I_{2n} 
            -h \mathcal H \left(\fiz I_{2n} -b\right) \\ 
       B_2 := \frac{\partial \Psi(z_0,z_h)}{\partial z_h} &=& - I_{2n} 
	    - h\mathcal H \left(\fiz I_{2n} + b\right) 
    \label{eqn:mats}
\end{eqnarray}
where $\mathcal H$ is a Hamiltonian matrix given by $\mathcal H= J D^2H$ and $D^2H$ is the Hessian matrix of $H$. 
The symplecticity test $ B_1JB_1^T - B_2JB_2^T = 0$ is satisfied since 
$\mathcal H$ and $b$ are Hamiltonian matrices.

\subsection{The generalized implicit symplectic Euler scheme}

The relevance of this method is that we have a linear (continuous) space of dimension 
$n(2n+1)$ where we can select Hamiltonian matrices for the construction of numerical schemes.
Using Theorem \ref{teo:main} we are able to define the \emph{generalized implicit symplectic 
Euler scheme} as the map given by 
\begin{eqnarray}
    \psi_h: U\times U & \to & U \nonumber\\
    (z_0, z_{h}) & \mapsto & z_{h} = z_0 + h X_H (\bar\bz)
    \label{eqn:first:appr}
\end{eqnarray}
where  $\bar\bz = \fiz(z_0+z_h) + b(z_h - z_0)$
and it corresponds to the map $\bar\bz=\rho(z_0,z_h)$ given by the expression (\ref{eqn:proj:coord:b}).

\begin{remark}
   In a slightly different context, Feng Kang (unpublished) has shown that generating functions 
   obtained by some particular type of matrices $\alpha\in M_{4n\times 4n}(\mathbb R)$ 
   can be reduced to generating functions constructed by an equivalent simplified matrix 
   containing the two submatrices given in expression 
   (\ref{eqn:proj:coord:b}) (see Ge and Dau-liu \cite{GD95}). In their method, those matrices
   are the input for the evolutive Hamilton-Jacobi equation. 
\end{remark}

\begin{proposition}
   The elementary Liouvillian form $\theta_0$ on an exact symplectic manifold $(M,\omega)$,
   produces a Liouvillian form on the product manifold 
   $(\tilde M,\omega_\ominus)$ given by $\theta = \pi_1^*\theta_0 - \pi_2^*\theta_0=\fiz dz_0 J z_0 - \fiz dz_h J z_h$. In 
   Darboux's coordinates it corresponds to
   \begin{eqnarray}
        \theta = \fiz\left(p_0dq_0 -q_0dp_0 -p_hdq_h + q_hdq_h\right),
   \end{eqnarray}
   whose generalized symplectic Euler scheme is the symplectic mid-point rule.
\end{proposition}
{\it Proof.} It is a direct computation using the fact that elementary Liouvillian form
does not have a symmetric component. In terms of the generalized symplectic Euler scheme,
their argument $\bar\bz$ have a null Hamiltonian matrix $b=0_{2n}$. 
The argument $\bar\bz = \fiz(z_0+z_h)$ gives the scheme
\begin{eqnarray}
    z_h = z_0 + h X_H\left( \fiz(z_0+z_h)\right),
\end{eqnarray}
which is the mid-point rule.
$\hfill\square$

\begin{proposition}
  In terms of the generalized symplectic Euler scheme (\ref{eqn:first:appr})
  with Darboux's coordinates $z_0=(q_0,p_0)^T$ and $z_h=(q_h,p_h)^T$, 
  the Hamiltonian matrices
\begin{eqnarray*}
    b_A = \fiz \left( 
     \begin{array}{c c}
       -I   & 0 \\
       0 &  I
    \end{array}
    \right),\qquad 
    b_B = \fiz \left( 
     \begin{array}{c c}
       I   & 0 \\
       0 &  -I
    \end{array}
    \right),\qquad b\in  \mathbb M_{2n\times 2n}(\mathbb R).
\end{eqnarray*} 
corresponds to the symplectic Euler $A$ and $B$ schemes, respectively.
\end{proposition}
{\it Proof.} Write the sum and difference vectors in Darboux's coordinates by 
\begin{eqnarray*}
    \fiz(z_0+z_h) = \fiz \left( 
     \begin{array}{c}
       q_0+q_h\\
       p_0+p_h
    \end{array}
    \right),\qquad {\rm and}\qquad
    (z_h-z_0) = \left( 
     \begin{array}{c}
       q_h-q_0\\
       p_h-p_0
    \end{array}
    \right).
\end{eqnarray*} 
Computing $b_A(z_h-z_0)$ leads to $b_A(z_h-z_0) = \fiz \left(q_0-q_h, p_h-p_0\right)^T$
and consequently $\bar \bz =\fiz(z_0+z_h) +b_A(z_h-z_0)$ produces $\bar \bz = (q_0, p_h)^T$.
In the same way $\bar \bz =\fiz(z_0+z_h) +b_B(z_h-z_0)$ produces $\bar\bz = (q_h,p_0)^T$.
Using these arguments in the generalized scheme (\ref{eqn:first:appr}) we obtain the 
standard Euler schemes.
$\hfill\square$

\section{Two families of implicit symplectic integrators}

In this section we develop two examples of the application of the 
method of Liouvillian forms for constructing symplectic integrators. 
The first one considers the construction of the special symplectic manifold 
on $\mathcal Q_1\times\mathcal Q_2$, using the symplectomorphism 
$\Psi:\tilde M\to T^*(\mathcal Q_1\times \mathcal Q_2)$.
The second one is given for applying the method directly. 
In the practice, we will use the generalized symplectic 
Euler scheme given in (\ref{eqn:first:appr}).

\subsection{Example 1: A simple family from symplectic rotations}
Here, we explicitely define the symplectomorphism between the product manifold 
$(\tilde M,\omega_\ominus)$ and the cotangent bundle $T^*(\mathcal Q_1\times \mathcal Q_2)$
using a family of symplectic rotations. We take the loop of symplectic rotations
induced by the complex structure $\tilde J\mapsto \exp(\phi \tilde J)$ on the product 
manifold $(\tilde M,\omega_\ominus)$. 
\begin{lemma}
   If $\tilde J$ is a complex structure then $\exp(\phi\tilde J)=\cos(\phi) I_{4n}+ \sin(\phi)\tilde J$.
\end{lemma}
{\it Proof.} Use $\tilde J^2=-I_{4n}$ in the development of the exponential. $\hfill\square$
Define the symplectic rotation $R_\phi \in Sp(\tilde M,\omega_\ominus)$ by $R_\phi = \exp(\phi\tilde J)$ 
and define the curve of diffeomorphisms $\Psi_\phi=E_1\circ R_\phi$,  $\phi\in[0,\pi/2]$.

Given symplectic coordinates $(q_0,p_0,q_h,p_h)\in\tilde M$ and $(x_0,x_h,y_0,y_h)\in T^*(\mathcal Q_1\times\mathcal Q_2)$,
the diffeomorphism $\Psi_{\phi}:\tilde{M}\to T^*(\mathcal Q_1\times\mathcal Q_2)$
pulls-back the Liouville form $\theta_{\mathcal Q_1 \times \mathcal Q_2}=y_0dx_0 + y_hdx_h$
into 
$\theta_{\phi} = \Psi^*_{\phi}\theta_{\mathcal Q_1\times\mathcal Q_2}$ expressed in these coordinates by
\begin{eqnarray}
    \theta_{\phi} &=& \cos^2\phi p_0 dq_0 - \sin^2\phi p_hdq_h - \sin^2\phi q_0dp_0 +\cos^2\phi q_hdp_h \nonumber \\
                & & + \cos\phi\sin\phi(q_0dq_0-p_0dp_0) + \cos\phi\sin\phi( q_hdq_h - p_h dp_h) 
                \label{eqn:theta:1}
\end{eqnarray}
The quintuple %we know %With these definitions 
$\left(\tilde{M}, (\mathcal Q_1\times\mathcal Q_2),  \theta_{\phi},\pi, \Psi_\phi \right)$ 
is a family of special symplectic manifolds on $\mathcal Q_1\times\mathcal Q_2$.
We have the following
\begin{lemma}
    The set of forms (\ref{eqn:theta:1})
    is a $1$-parameter family of Liouvillian forms. The path $\theta_{\phi}$ contains the Liouvillian forms 
    associated to the generating functions of type II and type III.
    \label{lem:11}
\end{lemma}
{\it Proof.}
By a direct computation we obtain $d\theta_\phi= dp_0\wedge dq_0 - dp_h\wedge dq_h = \omega_\ominus$,
showing that any element in the family is Liouvillian on $(\tilde M,\omega_\ominus)$.
The forms %associated to the generating functions of type \emph{II} and \emph{III} are 
at the values $\phi=0$ and $\phi=\frac{\pi}{2}$ are 
$$\theta_\phi|_{\phi=0} = p_0dq_0 + q_hdp_h\qquad{\rm and}\qquad
\theta_\phi |_{\phi=\pi/2} = -q_0dp_0 - p_hdq_h.$$ 
By Lemma \ref{lem:II:III} they are associated to 
generating functions of type \emph{II} and \emph{III}, respectively.
$\hfill\square$

The family of vertical fibers to $\Lambda_\phi\subset \tilde{M}$, 
associated to the Liouvillian forms $\theta_{\phi}$ is given by equations
\begin{eqnarray*}
    \hat q_0 = \cos^2\phi p_0 + \cos\phi\sin\phi q_0,\hspace{20pt}&\quad&
    \hat p_0 = - \sin^2\phi q_0 - \cos\phi_i\sin\phi p_0\\
    \hat q_h = - \sin^2\phi p_h + \cos\phi_i\sin\phi q_h,&\quad&
    \hat p_h = \cos^2\phi q_h - \cos\phi\sin\phi p_h.
\end{eqnarray*}
We obtain the tangent fibers by left multiplying vertical fibers by $\tilde J^T$,
mapping $(\hat q_0,\hat p_0, \hat q_h,\hat p_h)\mapsto (-\hat p_0,\hat q_0, \hat p_h,-\hat q_h)$. 
Finally, the projection $T\pi(T\Lambda_\phi)\cong T\pi_1(T\Lambda_\phi)\oplus  T\pi_2(T\Lambda_\phi)$
in local coordinates is $\rho(q_0,p_0,q_h,p_h)= \left(\hat p_h -\hat p_0, \hat q_0 -\hat q_h\right)$,  
or explicitely by
\begin{eqnarray}
    \bar{q_h} &=&  \cos^2\phi q_h + \sin^2\phi q_0 + \cos\phi\sin\phi( p_0 - p_h) \nonumber\\
    \bar{p_h} &=&  \cos^2\phi p_0 + \sin^2\phi p_h + \cos\phi\sin\phi( q_0 - q_h ).
    \label{eqn:proj}
\end{eqnarray}
The family of implicit symplectic integrators is given by 
\begin{eqnarray}
    \textstyle{ q_h = q_0 + h  \frac{\partial H}{\partial p}\left( \bar q_h, \bar p_h \right) } &\qquad& 
    \textstyle{ p_h = p_0 - h  \frac{\partial H}{\partial q}\left( \bar q_h, \bar p_h \right). }
\end{eqnarray}
Evaluating these coordinates in $\phi=0$ and $\phi=\pi/2$ we obtain 
\begin{eqnarray}
    \left.\left(\bar{q_h}, \bar{p_h}\right)\right|_{\phi=0} =  \left(q_h, p_0\right) \quad{\rm and}\quad
    \left.\left(\bar{q_h}, \bar{p_h}\right)\right|_{\phi=\pi/2} =  \left(q_0, p_h\right).
    \label{eqn:EulerAB}
\end{eqnarray}
We have proved the following result.
\begin{corollary}
  The implicit scheme $z_h = z_0 + h X_H\circ \rho(z_0,z_h)$
  where the linear map $\rho(z_0,z_h)=(\bar q_h, \bar p_h)$ is defined by the expression (\ref{eqn:proj}), is a family of 
  symplectic integrators which joins the symplectic Euler schemes A and B.
\end{corollary}
Using 
the identities 
$$\textstyle{\cos\phi\sin\phi = \frac{\sin2\phi}{2},\quad \sin^2\phi = \frac{1-\cos2\phi}{2},\quad
\cos^2\phi = \frac{1+\cos2\phi}{2}},$$
we rewrite (\ref{eqn:proj}) in compact form $\bar \bz = \left(\bar q_h, \bar p_h\right)$ by
$\bar \bz = \fiz(z_0+z_h)+b(z_h-z_0)$, where 
\begin{eqnarray}
    b = \fiz \left( 
     \begin{array}{c c}
       -\cos2\phi   & \sin2\phi \\
       \sin2\phi &  \cos 2\phi
    \end{array}
    \right),\qquad b\in  \mathbb M_{2n\times 2n}(\mathbb R).
    \label{eqn:b:2}
\end{eqnarray}

\subsection{Example 2: A three-parameter family}
The method of Liouvillian forms applied to the construction of symplectic
integrators yields the generalized symplectic Euler scheme (\ref{eqn:first:appr}).
This scheme uses the argument $\bar \bz = \fiz(z_0+z_h)+b(z_h-z_0)$, where 
$b\in \mathbb M_{2\times 2}(\mathbb R)$ is an arbitrary Hamiltonian matrix. 
In a problem with 1 degree of freedom, $b$ is a matrix of the form
\begin{eqnarray*}
    b = \left( 
     \begin{array}{c c}
       \alpha   & \beta \\
       \gamma &  -\alpha
    \end{array}
    \right),\qquad \alpha,\beta,\gamma\in\mathbb R.
\end{eqnarray*} 
These are three free parameters in the generalized Euler method.
Explicitely we have the numerical scheme 
\begin{eqnarray}
    q_h = q_0 + h\textstyle{\frac{\partial H}{\partial q}}\left( \bar q, \bar p \right) \qquad {\rm and} \qquad
    p_h = p_0 - h\textstyle{\frac{\partial H}{\partial p}}\left( \bar q, \bar p \right),
\end{eqnarray}
where 
\begin{eqnarray*}
  \bar q &=& \left(\fiz -\alpha\right)q_0 + \left(\fiz +\alpha\right)q_h +\beta(p_h-p_0)\\ 
  \bar p &=& \left(\fiz + \alpha\right)q_0 + \left(\fiz -\alpha\right)q_h + \gamma(q_h-q_0),
\end{eqnarray*}
and $\alpha,\beta,\gamma < h^2 < \fiz$ have small values. Parameter $\alpha$ modifies the symmetry of the oscillations
with respect to time and the limiting cases $(\alpha,\beta,\gamma)= (-\fiz,0,0)$ and  
$(\alpha,\beta,\gamma)=(\fiz,0,0)$ are the usual symplectic Euler $A$ and $B$ integrators. 

The parameters $\{\alpha,\beta,\gamma\}$ modify the error and the oscillations of the numerical
integrator around the orbits of the Hamiltonian flow. 
This fact has a geometrical justification that we develop in 
\cite{Jim16a}. In a parallel 
article \cite{JVR17}, we make a numerical study of this family explaining the 
oscillations of the numerical solution around the exact solution 
using a variational point of view.

% The Acknowledgements are an un-numbered section
\section*{Acknowledgements}
% Acknowledgements text here
The author thanks J.P. Vilotte and B. Romanowicz for their support and constructive
criticism on this work. 
This research was developed with support from the \emph{Fondation du Coll\`ege de 
France} and \emph{Total} under the research convention PU14150472, as well as the ERC Advanced Grant 
WAVETOMO, RCN 99285, Subpanel PE10 in the F7 framework.\\ 
Conflict of Interest: The author declares that he has no conflict of interest.


\begin{thebibliography}{10}
\providecommand{\url}[1]{{#1}}
\providecommand{\urlprefix}{URL }
\expandafter\ifx\csname urlstyle\endcsname\relax
  \providecommand{\doi}[1]{DOI~\discretionary{}{}{}#1}\else
  \providecommand{\doi}{DOI~\discretionary{}{}{}\begingroup
  \urlstyle{rm}\Url}\fi

 \bibitem{AM78}
Abraham, R., Marsden, J.: {Foundations of mechanics Second Ed.}
\newblock Benjamin Cummings (1978)

\bibitem{BGG03}
Bryant, R., Griffiths, P., Grossman, D.: Exterior differential systems and
  Euler-Lagrange partial differential equations.
\newblock University of Chicago Press (2003)

\bibitem{GD95}
Ge, Z., Dau-liu, W.: {On the invariance of generating functions for symplectic
  transformations}.
\newblock Diff. Geom. and its Appl \textbf{5}, 59--69 (1995)

\bibitem{Jim15f}
Jim{\'e}nez-P{\'e}rez, H., Vilotte, J.-P., Romanowicz, B.: 
{On the Poincar\'e's generating function and the
  symplectic mid-point rule}.
\newblock submitted arxiv:1508.07743v3  (2017).
\newblock \urlprefix\url{http://arxiv.org/abs/1508.07743v3}

\bibitem{Jim16a}
Jim{\'e}nez-P{\'e}rez, H.: {Hamilton-Liouville pairs}.
\newblock preprint  (2016)

\bibitem{JVR17}
Jim{\'e}nez-P{\'e}rez, H., Vilotte, J.-P., Romanowicz, B.: {The source of
  numerical oscillations in symplectic integration}.
\newblock submitted  (2017)

\bibitem{Lib00}
Libermann, P.: {On Liouville Forms}.
\newblock Poisson Geometry, Banach Center Publications \textbf{51}, 151--164
  (2000)

\bibitem{LM87}
Libermann, P., Marle, C.M.: {Symplectic Geometry and Analytical Mechanics}.
\newblock Ridel (1987)

\bibitem{MR91}
Marsden, J.E., Ratiu, T.S.: {Introduction to Mechanics and Symmetry}.
\newblock Springer-Verlag (1999)

\bibitem{Poi99}
Poincar{\'e}, H.: {Les m{\'e}thodes nouvelles de la m{\'e}canique c{\'e}leste
  Tome III}, vol. III.
\newblock Gauthier-Villars (1899)

\bibitem{STW92}
Simo, J., Tarnow, N., Wong, K.: {Exactly energy-momentum conserving algorithms
  and symplectic schemes for nonlinear dynamics}.
\newblock Comp. Meth. in Appl. Mech. and Engin \textbf{100}, 63--116 (1992)

\bibitem{ST72}
Sniatycki, J., Tulczyjew, W.: {Generating forms on Lagrangian submanifolds}.
\newblock Indiana Univ. Math. J \textbf{22} (1972)

\bibitem{Tul77}
Tulczyjew, W.: {The Legendre Transformation}.
\newblock Annales de l'IHP \textbf{section A}, 1 (101-114)

\bibitem{Tul76}
Tulczyjew, W.: {Les sous-vari{\'e}t{\'e}s lagrangiennes et la dynamique
  lagrangienne}.
\newblock C.R.Acad.Sci. Paris \textbf{283}, 675--678 (1976)

\bibitem{Wey46}
Weyl, H.: {The Classical Groups. Their Invariants and Representations.}
\newblock Princeton University Press (1946)
 

\end{thebibliography}
\end{document}